\documentclass[12pt,twoside]{elsart}

\usepackage{amsmath,amsfonts,amssymb,latexsym,texdraw}

\newcommand{\st}{${}\hspace{.5cm}$}

\begin{document}

\begin{frontmatter}

\title{Duality of translation association schemes coming from certain actions}
\author [dae]{Dae San Kim},
\author [hyun]{Hyun Kwang Kim}

\address[dae] {Department of Mathematics, Sogang University, Seoul 121-742, Korea}%
        \thanks{E-mail address : dskim@sogang.ac.kr (D.S.Kim)\\
        This work was supported by grant No. R01-2007-000-11176-0
        from the Basic Research Program of the Korea Science and Engineering Foundation.}
\address[hyun] {Department of Mathematics, Pohang University of Science and Technology, Pohang 790-784, Korea}

\begin{abstract}
\st Translation association schemes are constructed from actions of
finite groups on finite abelian groups satisfying certain natural
conditions. It is also shown that the mere existence of maps from
finite groups to themselves sending each element in their groups to
its \textquoteleft adjoint' entails the self-duality of the
constructed association schemes. Many examples of these, including
Hamming scheme and sesquilinear forms schemes, are provided. This
construction is further generalized to show the duality of the
association schemes coming from actions of two finite groups on the
same finite abelian group. An example of this is supplied with weak
Hamming schemes.
\end{abstract}
\begin{keyword}
Action ; Orbit ; Translation association scheme ; Adjoint ;$\\
{}\hspace{2.1cm}$Self-duality ; Duality
\end{keyword}
\end{frontmatter}

\section{Introduction}

\st Let $G$ be a finite group acting on a finite additive abelian
group $X$, and let $\mathcal{O}_{0}=\{0\}$, $\mathcal{O}_{1}, \cdots
, \mathcal{O}_{d}$ be the $G$-orbits. Assume that the action
satisfies the conditions
\begin{equation*}
  g(x+x')=gx+gx', ~for ~all ~g\in G ~and ~x, x'\in X,
\end{equation*}
and
\begin{equation*}
  x\in\mathcal{O}_{i} \Rightarrow -x\in\mathcal{O}_{i}, ~for ~0\leq i\leq
  d.
\end{equation*}
Then it is easy to see that
$\mathfrak{X}_{G}=(X,{~\{R_{i}\}}_{i=0}^{d})$, with $(x,y)\in
R_{i}\Leftrightarrow y-x\in \mathcal{O}_{i}$, is a translation
association scheme (cf. Theorem 2).\\
\st Assume now further that the following map sending each element
to its adjoint exists, i.e., there is a map $\iota : G \rightarrow
G$ such that
\begin{equation*}
  \langle gx,y\rangle = \langle x,\iota(g)y\rangle, ~for ~all ~g\in G,
  ~x,y\in X.
\end{equation*}
Here $\langle ~~, ~\rangle : X\times X\rightarrow
\mathbb{C}^{\times}$ is an inner product on $X$. Then it is shown
that such a simple requirement is enough to guarantee the
self-duality of $\mathfrak{X}_{G}$ (cf. Corollary 7). Many examples
of such pairs $(G,X)$ satisfying the above three conditions are
provided. For instance, the Hamming scheme and most of sesquilinear
forms
schemes fall within this category (cf. Section 4).\\
\st Suppose now that $\check{G}$ is another finite group acting on
the same finite abelian group $X$, with
$\check{\mathcal{O}_{0}}=\{0\}$, $\check{\mathcal{O}_{1}}, \cdots ,
\check{\mathcal{O}_{d}}$ the $\check{G}$-orbits, and that the
corresponding conditions to $(1)$ and $(2)$ for $\check{G}$ and
$\check{\mathcal{O}_{i}} ~(0\leq i\leq d)$ are satisfied. Assume in
addition that there is a map $\iota : G\rightarrow \check{G}$ such
that
\begin{equation*}
  \langle gx,y\rangle = \langle x,\iota(g)y\rangle , ~for ~all ~g\in G,
  ~x,y\in X.
\end{equation*}
Then it is shown that the existence of the map sending each element
to its adjoint is strong enough to yield the duality of
$\mathfrak{X}_G$ and $\mathfrak{X}_{\check{G}}$. An example of this
is illustrated with what we call the weak Hamming scheme $H(n_1,
\cdots , n_t, q)$. It is the wreath product of the Hamming schemes
$H(n_1,q), \cdots , H(n_t,q)$. Our result now says that $H(n_1,
\cdots , n_t, q)$ and $H(n_t, \cdots , n_1, q)$ are dual to each
other. This would have been observed already in \cite{HKKim} or even
in the earlier works \cite{DSKim:2007} and \cite{HKKim:2005}. The
weak Hamming scheme can be also constructed in connection with weak
order poset weight (cf. Section 6). Also, it is an example of weak
metric schemes which are the subject of the recent paper
\cite{DSKimGCKim}.\\

\section{Preliminaries}

\st A pair $\mathfrak{X}=(X,{~\{R_i\}}_{i=0}^d)$ consisting of a
finite set $X$ (referred to as the vertex set of $\mathfrak{X}$) and
$d+1$ nonempty subsets $R_i$ of $X\times X$ is called a $d$-class
(symmetric) association scheme if\\
{\rm(i)} {$\{R_0, R_1, \cdots , R_d \}$ is a partition of $X \times X$,}\\
{\rm(ii)} {$R_0 = \triangle_X = \{(x,x)|x\in X\}$,}\\
{\rm(iii)} {${}^tR_i = R_i$, ~for ~all ~$i$, where ${}^tR_i = \{(x,y)|(y,x)\in R_i\}$,}\\
{\rm(iv)} {for any $i,j,k(0\leq i,j,k \leq d)$, there are numbers,
            called intersection $
            {}\hspace{.8cm}$numbers, $p_{ij}^k$ such that, for any
            $(x,y)\in R_k$,
            \begin{equation*}
              p_{ij}^k = \# \{z\in X|(x,z)\in R_i, ~(z,y)\in R_j\}.
            \end{equation*}}\\
\st Let $\mathfrak{X} = (X,{~\{R_{i}\}}_{i=0}^{d})$ be an
association scheme. Let $A_i$ be the adjacency matrix of $R_i$, for
$0\leq i \leq d$. Then $A_0, A_1, \cdots , A_d$ generate a
$(d+1)$-dimensional commutative subalgebra $\mathcal{A}$ of
symmetric matrices in $M_n(\mathbb{C}) ~(n=|X|)$, called the
Bose-Mesner algebra of $\mathfrak{X}$. $\mathcal{A}$ has another
nice basis $E_0, E_1, \cdots , E_d$, called the irreducible
idempotents of $\mathcal{A}$. They are determined (up to permutation
of the indices $1, \cdots , d$) by :
\begin{align}
  \begin{array}{cc}
    \textrm{(i)} &{E_iE_j = \delta_{ij}E_i, ~\textrm{for all} ~i,j,}~~~~~~~~~~~~~~~~~~~~~~~~~~~~~~\\
    \textrm{(ii)} &{\sum_{i=0}^d E_i = I,}~~~~~~~~~~~~~~~~~~~~~~~~~~~~~~~~~~~~~~~~~~~~~~\\
    \textrm{(iii)} &{E_0=|X|^{-1}J ~(J ~\textrm{all-one matrix}),}~~~~~~~~~~~~~~~~~~~~\\
    \textrm{(iv)} &{E_0, E_1, \cdots , E_d ~\textrm{are linearly independent over}
            ~\mathbb{C}.}
  \end{array}
\end{align}
\st The $\mathbb{C}$-space generated by $A_0, A_1, \cdots , A_d$ is
also closed under the Hadamard multiplication $\circ$, with $J$ as
the multiplicative identity. Write
\begin{equation*}
  E_i \circ E_j = |X|^{-1}\sum_{k=0}^d q_{ij}^kE_k ~(0\leq i,j \leq
  d).
\end{equation*}
Then $q_{ij}^k$'s are actually nonnegative real numbers, called the
Krein parameters. Let
\begin{equation*}
  A_j = \sum_{i=0}^d p_{ij}E_i, ~~~E_j=|X|^{-1}\sum_{i=0}^d q_{ij}A_i.
\end{equation*}

The $p_{ij}$'s and $q_{ij}$'s are respectively called $p$- and
$q$-numbers. In particular, $v_i = p_{0i} = \#\{z|(x,z)\in R_i\}$,
for any fixed $x\in X$, and $m_i = q_{0i} =$ rank $E_i$ are
respectively called the valencies and the multiplicities. Also, $P =
(p_{ij})$ and $Q = (q_{ij})$ are respectively
called the first and the second eigenmatrix of $\mathfrak{X}$.\\
\st Let $X$ be a finite additive abelian group, and let
$\mathfrak{X} = (X,{~\{R_{i}\}}_{i=0}^{d})$ be a $d$-class
association scheme. Then $\mathfrak{X}$ is called a translation
association scheme if
\begin{equation*}
  (x,y)\in R_i \Rightarrow(x+z,y+z)\in R_i, ~for ~all ~z\in X ~and
  ~i.
\end{equation*}

Let
\begin{equation}
  X_i = \{x\in X|(0,x)\in R_i\}, ~for ~0\leq i \leq d.
\end{equation}

Then $X_0, X_1, \cdots , X_d$ give a partition of $X$, and
\begin{equation*}
  (x,y)\in R_i \Leftrightarrow y-x\in X_i ~(0\leq i \leq d).
\end{equation*}
\st The dual association scheme $\mathfrak{X}^* =
(X^*,{~\{R_{i}^*\}}_{i=0}^{d})$ of $\mathfrak{X}$ consists of the
group $X^*$ of characters on $X$ together with $d+1$ nonempty
subsets $R_i^*$ of $X^*\times X^*$ determined by :
\begin{equation*}
  (\chi,\psi)\in R_i^* \Leftrightarrow \psi\chi^{-1}\in X_i^*,
\end{equation*}
where $X_i^* = \{\chi\in X^*|E_i\chi = \chi\}$. Here $\chi$ is
viewed as the column vector with the $x$-component $(x\in X)$ given
by $\chi(x)$. For the proof of the following theorem, see [3,
Theorem 2.10.10].


\begin{thm}\label{thm1}
Let $\mathfrak{X} = (X,{~\{R_{i}\}}_{i=0}^{d})$ be a translation
association scheme with parameters
$p_{ij}^k,~q_{ij}^k,~v_i,~m_i,~P,~Q$.
Then we have the following :\\
{\rm(a)} {The dual scheme $\mathfrak{X}^*
            =(X,~{\{R_{i}^*\}}_{i=0}^{d})$
            is also a translation association$\\
            {}\hspace{.7cm}$scheme with parameters ${p^*}_{ij}^k = q_{ij}^k,
            ~{q^*}_{ij}^k = p_{ij}^k, ~v_i^* = m_i, ~m_i^* = v_i,\\
            {}\hspace{.7cm}~P^* = Q, ~Q^* = P$.}\\
{\rm(b)} {{\rm(i)} {$p_{ij} = \sum_{x\in X_j}\chi(x), ~for ~\chi\in X_i^*$,}$\\
            {}\hspace{.7cm}${\rm(ii)} {$q_{ij} = \sum_{\chi\in X_j^*}\chi(x), ~for ~x\in X_i$,}$\\
            {}\hspace{.7cm}${\rm(iii)} {$E_j = |X|^{-1}\sum_{\chi\in X_j^*}\chi~{}^t\overline{\chi}$,}$\\
            {}\hspace{.7cm}${\rm(iv)} {$v_j = |X_j|$,}$\\
            {}\hspace{.7cm}${\rm(v)} {$m_j = |X_j^*|$.}}\\
\end{thm}
\st Two association schemes $\mathfrak{X} =
(X,~{\{R_{X,i}\}}_{i=0}^{d})$ and $\mathfrak{Y} =
(Y,~{\{R_{Y,i}\}}_{i=0}^{d})$ are said to be isomorphic if there are
a bijection $f : X\rightarrow Y$ and a permutation $\sigma$ of $\{1,
2, \cdots , d\}$ such that
\begin{equation*}
  (x,y)\in R_{X,i} \Leftrightarrow (f(x),f(y))\in R_{Y,\sigma(i)},
  ~for ~1\leq i\leq d.
\end{equation*}
Here we always assume that $R_{X,0} = \triangle_X$, $R_{Y,0} =
\triangle_Y$, so that $(x,y)\in R_{X,0} \Leftrightarrow
(f(x),f(y))\in R_{Y,0}$. In particular, two $d$-class translation
association schemes $\mathfrak{X} = (X,~{\{R_{X,i}\}}_{i=0}^{d})$
and $\mathfrak{Y} = (Y,~{\{R_{Y,i}\}}_{i=0}^{d})$ are isomorphic if
there is a group isomorphism $f : X\rightarrow Y$ and a permutation
$\sigma$ of $\{1, 2, \cdots , d\}$ such that
\begin{equation*}
  x\in X_i \Leftrightarrow f(x)\in Y_{\sigma(i)}, ~for ~1\leq i\leq
  d.
\end{equation*}
Here $X_i = \{x\in X|(0,x)\in R_{X,i}\}$, $Y_i = \{y\in Y|(0,y)\in
R_{Y,i}\}$. If two association schemes are isomorphic, we may assume
that all the parameters of the schemes are the same.\\
\st For further facts about association schemes, one is referred to
\cite{Bannai} and \cite{Brouwer}.

\section{Construction of translation association schemes}

\st Let $G$ be a finite group acting on a finite additive abelian
group $X$, with $\mathcal{O}_{0}=\{0\}$, $\mathcal{O}_{1}, \cdots ,
\mathcal{O}_{d}$ the $G$-orbits. We assume that this action
satisfies the conditions
\begin{equation}
  g(x+x')=gx+gx', ~for ~all ~g\in G, ~x, x'\in X,
\end{equation}
and
\begin{equation}
  x\in\mathcal{O}_{i} \Rightarrow -x\in\mathcal{O}_{i}, ~for ~0\leq i\leq
  d.
\end{equation}
Note here that $g0 = 0$, for all $g\in G$, as
$\mathcal{O}_{0}=\{0\}$. This together with (3) implies that
\begin{equation}
  g(-x) = -gx, ~for ~all ~g\in G, ~x\in X.
\end{equation}\\
\st A special case of the following theorem appeared in
\cite{DSKim}.


\begin{thm}\label{thm2}
$\mathfrak{X}_G = (X,~{\{R_{i}\}}_{i=0}^{d})$, given by
\begin{equation*}
  (x,y)\in R_i\Leftrightarrow y-x\in \mathcal{O}_{i} ~(0\leq i\leq
  d),
\end{equation*}
is a translation association scheme.
\end{thm}
\begin{pf}
Here the conditions (3) and (4) are needed. All are easy to check,
perhaps except for the condition on the intersection numbers. Let
$x, y, x', y'\in X$, with $u = y-x, ~v = y'-x'\in \mathcal{O}_{k}$.
Then we must see :
\begin{equation*}
  \#\{z\in X|z-x\in \mathcal{O}_{i}, ~y-z\in \mathcal{O}_{j}\} =
  \#\{z\in X|z-x'\in \mathcal{O}_{i}, ~y'-z\in \mathcal{O}_{j}\}.
\end{equation*}
Observe that
\begin{equation*}
  \{z\in X|z-x\in \mathcal{O}_{i}, ~y-z\in
  \mathcal{O}_{j}\}\rightarrow \{z\in X|u-z\in \mathcal{O}_{i}, ~z\in
  \mathcal{O}_{j}\} ~(z\mapsto y-z)
\end{equation*}
is a bijection. So it is enough to show :
\begin{equation*}
  \#\{z\in X|u-z\in \mathcal{O}_{i}, ~z\in \mathcal{O}_{j}\} =
  \#\{z\in X|v-z\in \mathcal{O}_{i}, ~z\in \mathcal{O}_{j}\}.
\end{equation*}
As $G$ acts transitively on each $\mathcal{O}_{i}$, and $u,v\in
\mathcal{O}_{k}$, there is an $h\in G$ such that $hu = v$, and hence
\begin{equation*}
  \{z\in X|u-z\in \mathcal{O}_{i}, ~z\in
  \mathcal{O}_{j}\}\rightarrow \{z\in X|v-z\in \mathcal{O}_{i}, ~z\in
  \mathcal{O}_{j}\} ~(z\mapsto hz)
\end{equation*}
is a bijection. Notice that (3), (4) and (5) are used here.$\qed$
\end{pf}


\begin{rem}\label{remark1}
{\rm(1)} {For the association scheme $\mathfrak{X}_G = (X,~{\{R_{i}\}}_{i=0}^{d})$,
            \begin{equation*}
              X_i = \mathcal{O}_i, ~for ~0\leq i\leq d ~(cf. ~(2)).\end{equation*}}\\
{\rm(2)} {Let the condition (4) be replaced by :
            \begin{equation}
              for ~each ~i ~(0\leq i\leq d),
              ~x\in \mathcal{O}_i\Rightarrow -x\in \mathcal{O}_j, ~for ~some ~j.
            \end{equation}}
\end{rem}
If (3) is satisfied along with (6), then in fact we have $x\in
\mathcal{O}_i\Leftrightarrow -x\in \mathcal{O}_j$, and hence the
association scheme $\mathfrak{X} = (X, ~{\{R_{i}\}}_{i=0}^{d})$,
given by $(x,y)\in R_i\Leftrightarrow y-x\in \mathcal{O}_i ~(0\leq
i\leq d)$, is not a symmetric and yet commutative association
scheme. However, in this paper we will consider only the association
schemes coming from the actions satisfying (3) and
(4).\\
\par
\st Let $\langle ~~, ~\rangle : X\times X \rightarrow
\mathbb{C}^{\times}$ be an
inner product on the finite additive abelian group $X$, i.e.,\\
{\rm(i)} {$\langle x, y\rangle = \langle y, x\rangle$, for all $x,y\in X$,}\\
{\rm(ii)} {$\langle x,y+z\rangle = \langle x,y\rangle \langle x,z\rangle$, for all $x,y,z\in X$,}\\
{\rm(iii)} {$\langle x,y\rangle = \langle x,z\rangle$, for all $x\in X\Rightarrow y=z$.}\\


\begin{rem}\label{remark2}
It is  well-known that such an inner product always exists on a
finite abelian group $X$. $\langle ~~,x\rangle$ will denote the
character on $X$ given by $y\mapsto \langle y,x\rangle$, so that
\begin{equation*}
  X^* = \{\langle ~~,x\rangle|x\in X\}.
\end{equation*}

Also, $\langle ~~,x\rangle$ will indicate the column vector of size
$|X|$ whose $y$-component is $\langle y,x\rangle ~(y\in X)$.
\end{rem}


\begin{thm}\label{thm3}
The following are equivalent.\\ {\rm(a)} {$\mathfrak{X}_G =
(X,~{\{R_{i}\}}_{i=0}^{d})\rightarrow
            \mathfrak{X}_G^* = (X^*,~{\{R_{i}^*\}}_{i=0}^{d})$,
            given by $x\mapsto \langle ~~,x\rangle$, is an $\\
            {}\hspace{.7cm}$isomorphism of translation association schemes.}\\
{\rm(b)} {There is a permutation $\sigma$ of $\{1, 2, \cdots, d\}$ such that
            \begin{equation*}
              x\in X_j \Leftrightarrow \langle ~~,x\rangle\in
              X_{\sigma(j)}^*, ~for ~j=1, \cdots ,d.
            \end{equation*}}\\
{\rm(c)} {$E_j = |X|^{-1}\sum_{x\in X_j}\langle ~~,x\rangle
            ~{}^t\overline{\langle ~~,x\rangle} ~(0\leq j\leq d)$
            are the irreducible idempotents$\\
            {}\hspace{.7cm}$of the association scheme
            $\mathfrak{X}_G = (X,~{\{R_{i}\}}_{i=0}^{d})$.}\\
{\rm(d)} {$f_j : X\rightarrow \mathbb{C}$, given by $f_j(y) = \sum_{x\in X_j}\langle y,x\rangle$,
            is constant on each $X_i \\
            {}\hspace{.7cm}(0\leq i\leq d)$, for all $j=0, \cdots , d$.}\\
\end{thm}
\begin{pf}
$(a)\Leftrightarrow(b)$  This is just the definition.\\
$(a)\Rightarrow(d)$  As $(a)\Leftrightarrow(b)$, $q_{i\sigma(j)} =
\sum_{x\in X_j}\langle y,x\rangle$, for $y\in X_i$, is the
$q$-number of $\mathfrak{X}$, and hence is constant on each $X_i
~(0\leq i\leq d)$, for all $j=0, 1, \cdots , d$.\\
$(d)\Rightarrow(c)$  Let $\overline{E_j} = |X|^{-1}\sum_{x\in
X_j}\langle ~~,x\rangle ~{}^t\overline{\langle ~~,x\rangle}$, for
all $j$. Then $\overline{E_j} = |X|^{-1}\sum_{i=0}^d
\overline{q_{ij}}A_i$, with $\overline{q_{ij}} = \sum_{x\in
X_j}\langle y,x\rangle ~(y\in X_i)$. So $\overline{E_j}$ belongs to
the Bose-Mesner algebra of $\mathfrak{X}_G$.
\begin{equation*}
  (\overline{E}_i\overline{E}_j)_{ab} = |X|^{-2}\sum_{{}^{x\in X_i}_{y\in X_j}}
  \langle a,x\rangle\overline{\langle b,y\rangle}
  ~{}^t\overline{\langle ~~,x\rangle}\langle ~~,y\rangle.
\end{equation*}
Note that
\begin{equation*}
  {}^t\overline{\langle ~~,x\rangle}\langle ~~,y\rangle =
  \sum_{z\in X}\langle z,y-x\rangle =
    \begin{cases}
    ~0~,~&\text{}y\neq x,\\
    |X|,~&\text{}y = x.
    \end{cases}
\end{equation*}
So
\begin{equation*}
  (\overline{E}_i\overline{E}_j)_{ab} =
    \begin{cases}
    ~0~,\text{}i\neq j,\\
    |X|^{-1}\sum_{x\in X_i}\langle a-b,x\rangle =
        (\overline{E}_i)_{ab}, ~\text{}i = j,
    \end{cases}
\end{equation*}

and hence $\overline{E}_i\overline{E}_j =
\delta_{ij}\overline{E}_i$. Similarly, $\sum_{i=0}^d\overline{E}_i =
I$. Assume that $\sum_{i=0}^d\alpha_i\overline{E}_i = 0$, with
$\alpha_i\in \mathbb{C}$. Then $\alpha_i\overline{E}_i = 0$, for all
$i$. To show independence of $\overline{E}_0, \cdots ,
\overline{E}_d$, it is enough to see that $\overline{E}_i \neq 0$,
for all $i$. This is indeed the case, since
\begin{equation}
  \overline{E}_i\langle ~~,x\rangle =
    \begin{cases}
    \langle ~~,x\rangle,~&\text{if}~x\in X_i,\\
    ~~0~, ~&\text{if}~x\notin X_i.
    \end{cases}
\end{equation}
Thus $E_0, E_1, \cdots , E_d$ are the irreducible idempotents of
$\mathfrak{X}$ (cf. (1)). Note here
that $\overline{E}_0 = |X|^{-1}J$.\\
$(c)\Rightarrow(b)$  Recall that $\widetilde{E}_0, \widetilde{E}_1,
\cdots , \widetilde{E}_d$, with $\widetilde{E}_j =
|X|^{-1}\sum_{\chi\in X_j^*}\chi ~{}^t\overline{\chi}$, is also the
irreducible idempotents of $\mathfrak{X}_G$ (cf. Theorem 1, (b)).
Observe also that $\widetilde{E}_0 = |X|^{-1}J$. As the irreducible
idempotents are unique up to permutation, there is a permutation
$\sigma$ of $\{1, 2, \cdots ,d\}$ such that $E_j =
\widetilde{E}_{\sigma(j)}$, for $j = 1, \cdots ,d$. Now, using
(7) we have :\\
\begin{align*}
  x\in X_j &\Leftrightarrow E_j \langle ~~,x\rangle = \langle
  ~,x\rangle\\
        &\Leftrightarrow \widetilde{E}_{\sigma(j)}\langle
  ~,x\rangle = \langle ~~,x\rangle\\
        &\Leftrightarrow \langle ~~,x\rangle\in X_{\sigma(j)}^*.\qed
\end{align*}
\end{pf}

\st Assume now further that there is a map $\iota : G\rightarrow G$
such that
\begin{equation}
  \langle gx,y\rangle = \langle x,\iota(g)y\rangle, ~for ~all ~g\in G,
  ~x,y\in X,
\end{equation}
where $\langle ~~,~\rangle : X\times X \rightarrow
\mathbb{C}^\times$ is an inner product.


\begin{lem}\label{lemma1}
Under the assumption of (8), the sum
\begin{equation*}
  \sum_{x\in X_j}\langle y,x\rangle ~(y\in X_i)
\end{equation*}
depends only on $i$, for all $i,j$ with $0\leq i,j  \leq d$.
\end{lem}

\begin{pf}
Let $y_1, y_2\in X_i$. Then $y_2 = hy_1$, for some $h\in G$. So
\begin{align*}
  \sum_{x\in X_j}\langle y_2,x\rangle &= \sum_{x\in X_j}\langle hy_1,x\rangle\\
        &= \sum_{x\in X_j}\langle y_1,\iota(h)x\rangle\\
        &= \sum_{x\in X_j}\langle y_1,x\rangle.\qed
\end{align*}
Now, we get the following corollary from Theorem 5 which says in
particular that $\mathfrak{X}_G$ is self-dual.
\end{pf}


\begin{cor}\label{corollary1}
Let $\mathfrak{X}_{G} = (X,{~\{R_{i}\}}_{i=0}^{d})$ be the
translation association scheme obtained from the action of the
finite group $G$ on the finite abelian group $X$, satisfying (3) and
(4). Assume that there is a map $\iota : G\rightarrow G$ such that
\begin{equation*}
  \langle gx,y\rangle = \langle x,\iota(g)y\rangle, ~for ~all ~g\in G,
  ~x,y\in X,
\end{equation*}
where $\langle ~~,~\rangle : X\times X \rightarrow
\mathbb{C}^{\times}$ is an inner product. Then\\
{\rm(a)} {$\mathfrak{X}_{G} = (X,{~\{R_{i}\}}_{i=0}^{d})\rightarrow
            \mathfrak{X}_{G}^* = (X^*,{~\{R_{i}^*\}}_{i=0}^{d})$,
            given by $x\mapsto \langle ~~,x\rangle$, is an$\\
            {}\hspace{.7cm}$ isomorphism, i.e., $\mathfrak{X}_G$ is self-dual,}\\
{\rm(b)} {$E_j = |X|^{-1}\sum_{x\in X_j}\langle
            ~~,x\rangle~{}^t\overline{\langle ~~,x\rangle}
            ~(0\leq j\leq d)$ are the irreducible idempotents$\\
            {}\hspace{.7cm}$for $\mathfrak{X}_{G}$,}\\
{\rm(c)} {$q_{ij} = \sum_{x\in X_j}\langle y,x\rangle ~(y\in X_i)$
            are the $q$-numbers for $\mathfrak{X}_{G}$,}\\
{\rm(d)} {$p_{ij}^k = q_{ij}^k, ~p_{ij} = q_{ij}, ~v_i = m_i$.}\\
\end{cor}

\section{Examples for Section 3}

\st Here we will demonstrate that there are abundant examples of
actions of finite groups $G$ on finite abelian groups $X$ satisfying
(3) and (4), and (8) for suitable inner products on $X$. So the
schemes $\mathfrak{X}_{G} = (X,{~\{R_{i}\}}_{i=0}^{d})$ constructed
from these actions are, in particular, self-dual. In below,
$\lambda$ will always denote a fixed nontrivial additive character
on $\mathbb{F}_q$. The
examples (a) and (b) below are adopted from \cite{Delsarte:1973}.\\
\par
{\rm(a)} {Let $X$ be a finite abelian group with period $\nu$.
            Then $(\mathbb{Z}/(\nu))^{\times}$ acts on $X$ via
            \begin{equation*}
              (\mathbb{Z}/(\nu))^{\times} \times X\rightarrow X
              ~((\overline{m},x)\mapsto mx).
            \end{equation*}

            This action satisfies (3) and (4), and the orbits $\mathcal{O}_0 = \{0\}$,
             $\mathcal{O}_1, \cdots , \mathcal{O}_d$ are called the central classes of $X$.
            Further, given any inner product on $X$, (8) is satisfied with $\iota$ the identity map.}\\
\par
{\rm(b)} {Let $\omega$ be a primitive element in $\mathbb{F}_q$
            ($q$ an odd prime power), and let $d$ be a positive
            integer such that $2d\mid q-1$. Let $G = \langle
            \omega^d\rangle$ be the cyclic subgroup of
            $\mathbb{F}_q^{\times}$ of order $r = (q-1)/d$. Then
            $G$ acts on $X = (\mathbb{F}_q,+)$ by left
            multiplication. Here the orbits are $\mathcal{O}_0 = \{0\}$,
             $\mathcal{O}_1, \cdots , \mathcal{O}_d$, where, for
            $0\leq i\leq d-1$,
            \begin{equation*}
              \mathcal{O}_{i+1} = \{\omega^i, \omega^{d+i}, \cdots ,
              \omega^{(r-1)d+i}\}.
            \end{equation*}
            $\mathcal{O}_0, \mathcal{O}_1, \cdots , \mathcal{O}_d$ are called the cyclotomic
            classes of $\mathbb{F}_q$. The condition (3) is
            obviously satisfied, and the condition (4) is also
            valid, as we assume $2d\mid q-1$. Then $\langle x,y\rangle
             = \lambda(xy)$ is an inner product on $X$, and (8) is
            satisfied with $\iota$ the identity map.}\\
\par
\st The examples (c)-(f) will be about sesquilinear forms
association schemes. There are many articles about this topic. Here
we are content with just mentioning [3, Sections 9.5-6] and
[5,7,8,16].\\
\par
{\rm(c)} {Let $X = (\mathbb{F}_q^{m\times n},+)$, with $m\leq n$,
and let $G = GL(m,q)\times GL(n,q)$
            be the direct product of general linear groups. $G$ now acts on $X$ via
            \begin{equation*}
              G\times X \rightarrow X ~(((\alpha, \beta),A)\mapsto \phi^{(\alpha,\beta)}A := {}^t\alpha
              A\beta).
            \end{equation*}
            Then $\mathcal{O}_i = \{A\in X|rank(A) = i\} ~(i = 0, 1, \cdots ,
            m)$ are the $G$-orbits, and the conditions (3) and (4)
            hold. The associated scheme $\mathfrak{X}_{G}$ is
            called the bilinear forms scheme, which is usually
            denoted by $Bil(m\times n,q)$. Moreover, for $A\cdot B
            = tr(A~{}^tB) = \sum_{i,j}A_{ij}B_{ij} ~(A,B\in X)$, and
            $(\alpha,\beta)\in G$,
            \begin{equation*}
              \phi^{(\alpha, \beta)}A\cdot B = A\cdot
              \phi^{({}^t\alpha,{}^t\beta)}B.
            \end{equation*}
            So, for the inner product $\langle A,B\rangle = \lambda(A\cdot
            B)$ and $\iota : G\rightarrow G ~((\alpha, \beta)\mapsto
            ({}^t\alpha,{}^t\beta))$, (8) is valid.
            }\\
\par
{\rm(d)} {Let $X$ be the group of all alternating matrices of order
$m$
            over $\mathbb{F}q$. Recall here that $(A_{ij})$ is
            alternating $\Leftrightarrow A_{ii} = 0$, for $1\leq
            i\leq m$, and $A_{ji} = -A_{ij}$, for $1\leq i < j\leq
            m$. $G = GL(m,q)$ acts on $X$ via
            \begin{equation*}
              G\times X\rightarrow X ~((\alpha,A)\mapsto \phi^{(\alpha)}A := {}^t\alpha
              A\alpha).
            \end{equation*}
            Then $\mathcal{O}_i = \{A\in X|rank(A) = 2i\} ~(i = 0, 1, \cdots ,
            n=\lfloor \frac{m}{2}\rfloor)$ are the $G$-orbits, and
            the conditions (3) and (4) are satisfied. The
            associated scheme $\mathfrak{X}_{G}$ is called the
            alternating forms scheme which is denoted by
            $Alt(m,q)$. For $A\cdot B = \sum_{i<j}A_{ij}B_{ij} ~(A,B\in
            X)$, and $\alpha\in G$,
            \begin{equation*}
              \phi^{(\alpha)}A\cdot B = A\cdot
              \phi^{({}^t\alpha)}B.
            \end{equation*}
            This holds regardless of the characteristic of
            $\mathbb{F}_q$. But it is easier to check this for
            char $\mathbb{F}_q \neq 2$, since $A\cdot B =
            2^{-1}trA~{}^tB$ in that case. So, for the inner
            product $\langle A,B\rangle = \lambda(A\cdot B)$, and
            $\iota : G\rightarrow G ~(\alpha\mapsto {}^t\alpha)$, (8)
            is satisfied.
            }\\
\par
{\rm(e)} {Let $X$ be the group of all Hermitian matrices of order
$m$ over $\mathbb{F}_{q^2}$.
            Recall here that $A\in X$ is Hermitian if ${}^*A = A$,
            with $({}^*A)_{ij} = \overline{A_{ji}} = (A_{ji})^q$.
            $G = GL(m,q^2)$ acts on $X$ via
            \begin{equation*}
              G\times X\rightarrow X ~((\alpha,A)\mapsto
              \phi^{(\alpha)}A = {}^*\alpha A\alpha).
            \end{equation*}
            Then $\mathcal{O}_i = \{A\in X|rank(A) = i\} ~(i = 0, 1, \cdots ,
            m)$ are the $G$-orbits, and the conditions (3) and (4) are
            satisfied. The associatied scheme $\mathfrak{X}_{G}$ is called
            the Hermitian forms scheme which is denoted by
            $Her(m,q^2)$. For $A\cdot B =
            \sum_{i,j}A_{ij}\overline{B_{ij}} = tr(A~{}^*B) = tr(AB) ~(A,B\in
            X)$, and $\alpha\in G$,
            \begin{equation*}
              \phi^{(\alpha)}A\cdot B = A\cdot \phi^{({}^*\alpha)}B.
            \end{equation*}
            Note here that $A\cdot B\in \mathbb{F}_q$, and hence
            that $\langle A,B\rangle = \lambda(A\cdot B)$ makes
            sense and $\langle A,B\rangle$ is an inner product on
            $X$. Now, (8) holds for $\iota : G\rightarrow G
            ~(\alpha\mapsto {}^*\alpha)$.
            }\\
\par
{\rm(f)} {Let $X$ be the group of all symmetric matrices of order
$m$ over $\mathbb{F}_q$. Let $q$ be odd.
            $G = GL(m,q)$ acts on $X$ via
            \begin{equation*}
              G\times X\rightarrow X ~((\alpha,A)\mapsto
              \phi^{(\alpha)}A := {}^t\alpha A\alpha).
            \end{equation*}
            Then $\mathcal{O}_0 = \{0\}$, $\mathcal{O}_{r,+} = \{A\in X|A\sim J_r^+\}$,
            $\mathcal{O}_{r,-} = \{A\in X|A\sim J_r^-\} ~(r = 1, \cdots ,
            m)$ are the $G$-orbits, where $J_r^+ = I_r\dotplus O$, $J_r^- = \varepsilon I_1\dotplus I_{r-1}\dotplus
            O$, with $\varepsilon$ a fixed nonsquare element in
            $\mathbb{F}_q$. Here \textquoteleft$\sim$' and \textquoteleft$\dotplus$'
            indicate respectively \textquoteleft cogredient' and the matrix
            direct sum. Also, for this well-known fact one is
            referred to [15, Chap.IV]. The condition (3) is clearly satisfied.
            Assume further that $q \equiv 1 \pmod{4}$. Then, as $-1$ is a square,
            (4) is also valid. The associated scheme $\mathfrak{X}_{G}$ is called
            the symmetric forms scheme. For $A\cdot B =
            \sum_{i,j}A_{ij}B_{ij} = tr(A~{}^tB) = tr(AB) ~(A,B\in
            X)$, and $\alpha\in G$,
            \begin{equation*}
              \phi^{(\alpha)}A\cdot B = A\cdot \phi^{({}^t\alpha)}B.
            \end{equation*}
            Observe here that $\langle A,B\rangle = \lambda(A\cdot B)$
            is indeed an inner product on $X$, since we assume $q$
            is odd. Now, (8) holds for $\iota : G\rightarrow G
            ~(\alpha\mapsto {}^t\alpha)$. On the other hand, if
            $q\equiv 3 \pmod{4}$, then (4) is not satisfied but
            (6) is. So in that case the associated scheme is a
            commutative association scheme (cf. Remarks 3).
            }\\
\par
{\rm(g)} {Let $X = (\mathbb{F}_q^n,+)$, with $w_H$ the Hamming
            weight on $X$.
            Let $G = Aut(X,w_H)$ be the subgroup of $Aut(X)$
            consisting of all linear automorphisms $\phi$ of $X$
            preserving the Hamming weight, i.e., $w_H(\phi u) =
            w_H(u)$, for all $u\in X$. Then, as is well-known, $G
            \cong S_n ~{}_\psi\ltimes(\mathbb{F}_q^{\times})^n$,
            where $\psi : S_n\rightarrow
            Aut((\mathbb{F}_q^{\times})^n)$ is given by
            \begin{equation*}
              \sigma\mapsto (\alpha = (\alpha_1, \cdots ,
              \alpha_n)\mapsto \alpha_{\sigma^{-1}} = (\alpha_{\sigma^{-1}(1)}, \cdots ,
              \alpha_{\sigma^{-1}(n)})).
            \end{equation*}
            The isomorphism $S_n{}_\psi\ltimes(\mathbb{F}_q^{\times})^n
            \widetilde{\rightarrow} G$ is given by
            \begin{equation*}
              (\sigma,\alpha = (\alpha_1, \cdots ,
              \alpha_n))\mapsto \phi_{\sigma}\phi_{\alpha},
            \end{equation*}
            where
            \begin{equation*}
              \phi_{\sigma}(x_1, \cdots , x_n) =
              (x_{\sigma^{-1}1}, \cdots , x_{\sigma^{-1}n}), ~\phi_{\alpha}(x_1, \cdots , x_n) =
              (\alpha_1 x_1, \cdots , \alpha_n x_n).
            \end{equation*}
            Now, $G$
            acts naturally on $X$, and the $G$-orbits are
            $\mathcal{O}_i = \{x\in X|w_H(x) = i\} ~(i = 0, 1, \cdots ,
            n)$. Clearly, (3) and (4) are satisfied. The
            associated scheme $\mathfrak{X}_{G}$ is nothing but
            the Hamming scheme $H(n,q)$. For the usual
            $\mathbb{F}_q$-valued inner product $x\cdot y =
            \sum_{i=1}^{n}x_iy_i$ on $\mathbb{F}_q^n$, we observe
            \begin{equation*}
              \phi_\sigma \phi_\alpha x\cdot y =
              \sum_{i}\alpha_{\sigma^{-1}i}x_{\sigma^{-1}i}y_i =
              \sum_{i}x_i \alpha_i y_{\sigma i} = x\cdot
              \phi_{\sigma^{-1}}\phi_{\alpha_{\sigma^{-1}}}y.
            \end{equation*}
            So, for the inner product $\langle x,y\rangle =
            \lambda(x\cdot y)$, and $\iota : G\rightarrow G ~(\phi_\sigma \phi_\alpha\mapsto
            \phi_{\sigma^{-1}}\phi_{\alpha_{\sigma^{-1}}})$, (8) is
            valid.
            }\\

\section{Further generalizations}

\st The construction in Section 3 will be further generalized. Let
$X$ be a finite additive abelian group, and let $G$, $\check{G}$ be
two finite groups acting on $X$ with their respective orbits
$\mathcal{O}_0 = \{0\}$, $\mathcal{O}_1, \cdots , \mathcal{O}_d$ and
$\check{\mathcal{O}}_0 = \{0\}$, $\check{\mathcal{O}}_1, \cdots ,
\check{\mathcal{O}}_d$. Assume both actions satisfy the conditions
(3) and (4), with $\mathfrak{X}_G = (X, ~{\{R_{i}\}}_{i=0}^{d})$,
$\mathfrak{X}_{\check{G}} = (X, ~{\{\check{R}_{i}\}}_{i=0}^{d})$
their respective associated schemes. As in (2), we let, for $0\leq
i\leq d$,
\begin{equation*}
  X_i = \{x\in X|(0,x)\in R_i\}, ~\check{X}_i = \{x\in X|(0,x)\in
  \check{R}_i\},
\end{equation*}
so that
\begin{equation*}
  X_i = \mathcal{O}_i, ~\check{X}_i = \check{\mathcal{O}}_i, ~\textrm{for
  ~all} ~0\leq i\leq d.
\end{equation*}\\
\st A special case of the following theorem appeared in
\cite{DSKim}. The proof is left to the reader, as it follows by
slightly modifying that of Theorem 5. In the next theorem, $\langle
~~,~\rangle : X\times X \rightarrow \mathbb{C}^{\times}$ is
a fixed inner product.\\


\begin{thm}\label{thm4}
The following are equivalent.\\
{\rm(a)} {$\mathfrak{X}_{\check{G}} = (X,
            ~{\{\check{R}_{i}\}}_{i=0}^{d}) \rightarrow
            \mathfrak{X}_{G}^* = (X^*, ~{\{{R}_{i}^*\}}_{i=0}^{d})$, given by $x\mapsto \langle
            ~~,x\rangle$, is an$\\
            {}\hspace{.7cm}$isomorphism of translation
            association schemes.
            }\\
{\rm(b)} {There is a permutation $\sigma$ of $\{1, 2, \cdots , d\}$
such that
            \begin{equation*}
              x\in \check{X}_j \Leftrightarrow \langle ~~,x\rangle\in
              X_{\sigma(j)}^*, ~for ~j = 1, \cdots , d.
            \end{equation*}
            }\\
{\rm(c)} {$E_j = |X|^{-1}\sum_{x\in \check{X}_j}\langle ~~,x\rangle
{}^t\overline{\langle ~~,x\rangle}
            ~(0\leq j\leq d)$ are the irreducible idempotents$\\
            {}\hspace{.7cm}$of the association
            scheme $\mathfrak{X}_G = (X, ~{\{R_{i}\}}_{i=0}^{d})$.
            }\\
{\rm(d)} {$f_j : X \rightarrow \mathbb{C}$, given by $f_j(y) =
            \sum_{x\in \check{X}_j}\langle y,x\rangle$ is constant
            on each $X_i ~(0\leq i\leq\\
            {}\hspace{.7cm} d)$, for all $j = 0, 1, \cdots , d$.\
            }\\
\end{thm}

\st Assume now further that there is a map $\iota : G\rightarrow
\check{G}$ such that
\begin{equation}
  \langle gx,y\rangle = \langle x,\iota(g)y\rangle, ~for ~all ~g\in G,
  ~x,y\in X.
\end{equation}
Then, as in Lemma 6, we have the following.


\begin{lem}\label{lemma2}
Under the assumption of (9), the sum
\begin{equation*}
  \sum_{x\in \check{X}_j}\langle y,x\rangle ~(y\in X_i)
\end{equation*}
depends only on $i$, for all $i,j$ with $0\leq i,j \leq d$.
\end{lem}

\st So we obtain the following corollary from Theorem 8 and Theorem
1.


\begin{cor}\label{corollary2}
Let $\mathfrak{X}_G = (X, ~{\{R_{i}\}}_{i=0}^{d})$ and
$\mathfrak{X}_{\check{G}} = (X, ~{\{\check{R}_{i}\}}_{i=0}^{d})$ be
the two translation association schemes obtained from the actions of
the finite groups $G$ and $\check{G}$ on the same finite abelian
group $X$, and satisfying (3) and (4). Assume that there is a map
$\iota : G\rightarrow \check{G}$ such that
\begin{equation*}
  \langle gx,y\rangle = \langle x,\iota(g)y\rangle, ~for ~all ~g\in G,
  ~x,y\in X,
\end{equation*}
where $\langle ~~,~\rangle : X\times X \rightarrow
\mathbb{C}^{\times}$ is an inner product. Then\\
{\rm(a)} {$\mathfrak{X}_{\check{G}} = (X,
            ~{\{\check{R}_{i}\}}_{i=0}^{d})\rightarrow \mathfrak{X}_G^* = (X^*, ~{\{R_{i}^*\}}_{i=0}^{d})
            ~(x\mapsto \langle ~~,x\rangle)$ is an isomorphism,$\\
            {}\hspace{.7cm}$i.e.,
            $\mathfrak{X}_G$ and $\mathfrak{X}_{\check{G}}$ are dual
            to each other,
            }\\
{\rm(b)} {$E_j = |X|^{-1}\sum_{x\in \check{X}_j}\langle ~~,x\rangle
            ~{}^t\overline{\langle ~~,x\rangle}
            ~(0\leq j\leq d)$ are the irreducible idempotents$\\
            {}\hspace{.7cm}$for
            $\mathfrak{X}_G$,
            }\\
{\rm(c)} {$q_{ij} = \sum_{x\in \check{X_j}}\langle y,x\rangle ~(y\in
            X_i)$ are the $q$-numbers for $\mathfrak{X}_G$,
            }\\
{\rm(d)} {$p_{ij} = \sum_{x\in X_j}\langle y,x\rangle ~(y\in
            \check{X}_i)$ are the $p$-numbers for $\mathfrak{X}_G$,
            }\\
{\rm(e)} {$p_{ij}^k = \check{q}_{ij}^k, ~p_{ij} = \check{q}_{ij},
            ~m_i = \check{v}_i, ~q_{ij}^k = \check{p}_{ij}^k,
            ~~q_{ij} = \check{p}_{ij}, ~v_i = \check{m}_i$.
            }\\
\end{cor}

\section{An example for Section 5}

\st Here we are content with giving only one example for Section 5
which is what we call the weak Hamming scheme. Let $H(m,q)$ denote,
as usual, the Hamming scheme whose vertex set is $\mathbb{F}_q^m$
and $i$-th relation is given by
\begin{equation*}
  (x,y)\in R_i \Leftrightarrow d_H(x,y) = w_H(x-y) = i ~( i = 0, 1, \cdots ,
  m).
\end{equation*}
Here $w_H$ and $d_H$ are respectively the Hamming weight and the
Hamming metric. Then the weak Hamming scheme $H(n_1, \cdots , n_t,
q)$ is given as the wreath product (cf. \cite{Bang})
\begin{equation*}
  H(n_1, \cdots , n_t, q) = H(n_1,q)\wr \cdots \wr H(n_t,q),
\end{equation*}
so that the vertex set of that is $\mathbb{F}_q^{n_1}\times \cdots
\times \mathbb{F}_q^{n_t}$, and
\begin{align*}
  (x,y)\in &R_{n_1+ \cdots +n_{i-1}+i_0}\\
  &\Leftrightarrow x_{i+1} =
  y_{i+1}, \cdots , x_t = y_t ~~and ~~d_H(x_i,y_i) = w_H(x_i-y_i) =
  i_0,
\end{align*}
for $i = 1, \cdots , t, ~1\leq i_0\leq n_i$, or $i = 1, ~i_0 = 0$.\\
\st There is another description of $H(n_1, \cdots , n_t, q)$ which
has to do with poset-weight (poset-metric). This notion of
poset-weight (poset-metric) was first introduced in \cite{Brualdi}.
Let $\mathbb{P} = ([n],\leq)$ be a poset, with $[n] = \{1, 2, \cdots
, n\}$. The $\mathbb{P}$-weight $w_{\mathbb{P}}$ is the function on
$\mathbb{F}_q^n$, which is given by :
\begin{equation*}
  w_{\mathbb{P}}(x) = \#\{i\in [n]|i\leq j, ~\textrm{for ~some} ~j\in
  Supp(x)\}.
\end{equation*}
Here $Supp(x) = \{j\in [n]|x_j \neq 0\}$, for $x = (x_1, \cdots ,
x_n)\in \mathbb{F}_q^n$. Then $d_{\mathbb{P}}(x,y) =
w_{\mathbb{P}}(x-y)$ is a distance function on $\mathbb{F}_q^n$,
called $\mathbb{P}$-metric. We now specialize $\mathbb{P}$ as the
weak order poset $\mathbb{P}_0$. $\mathbb{P}_0 = n_1
\mathbf{1}\bigoplus \cdots \bigoplus n_t \mathbf{1}$ is given as the
ordinal sum of the antichains $n_i \mathbf{1}$ on the set $\{n_1+
\cdots + n_{i-1}+1, \cdots , n_1+ \cdots + n_{i-1}+n_i\}$, for $i =
1, \cdots , t$, i.e., the underlying set is $[n] ~(n = n_1+ \cdots
+n_t)$ and the order relation is given by :
\begin{equation*}
  k<l\Leftrightarrow k\in n_i \mathbf{1}, ~l\in n_j \mathbf{1}, ~\textrm{for ~some} ~i<j.
\end{equation*}
Let $G = Aut(\mathbb{F}_q^n,w_{\mathbb{P}_0})$ be the group of all
linear automorphisms of $\mathbb{F}_q^n$ preserving
$w_{\mathbb{P}_0}$-weight, i.e., $w_{\mathbb{P}_0}(\phi u) =
w_{\mathbb{P}_0}(u)$, for all $u\in \mathbb{F}_q^n$. Then $G$ acts
on $\mathbb{F}_q^n$ in a natural way and the $G$-orbits are
$\mathcal{O}_i = \{x\in \mathbb{F}_q^n|w_{\mathbb{P}_0}(x) = i\}$,
$i = 0, 1, \cdots , n$. Clearly, the conditions (3) and (4) are
valid. Moreover, for $i = 1, \cdots , t, ~1\leq i_0\leq n_i$, or $i
= 1, ~i_0 = 0$,
\begin{align*}
  y-x\in \mathcal{O}_{n_1+ \cdots +n_{i-1}+i_0} &\Leftrightarrow
  w_{\mathbb{P}_0}(y-x) = n_1+ \cdots +n_{i-1}+i_0\\
    &\Leftrightarrow x_{i+1} = y_{i+1}, \cdots , x_t = y_t ~and
    ~d_H(x_i,y_i) = i_0\\
    &\Leftrightarrow (x,y)\in R_{n_1+ \cdots +n_{i-1}+i_0}.
\end{align*}
Here we identified $\mathbb{F}_q^n$ with $\mathbb{F}_q^{n_1}\times
\cdots \times \mathbb{F}_q^{n_t}$ by writing the elements $x\in
\mathbb{F}_q^n$ as the blocks of coordinates $x = (x_1, \cdots ,
x_t)\in \mathbb{F}_q^{n_1}\times \cdots \times \mathbb{F}_q^{n_t}$.
So the associated scheme $\mathfrak{X}_G$ is nothing but the weak
Hamming scheme $H(n_1, \cdots , n_t, q)$. Similarly, the weak
Hamming scheme $H(n_t, \cdots , n_1, q) = H(n_t, q)\wr H(n_{t-1},
q)\wr \cdots \wr H(n_1, q)$ is also obtained from the action of the
group $\check{G} = Aut(\mathbb{F}_q^n,w_{\check{\mathbb{P}}_0})$ on
$\mathbb{F}_q^n$. Here $\check{\mathbb{P}}_0 = n_t
\mathbf{1}\bigoplus n_{t-1} \mathbf{1}\bigoplus \cdots \bigoplus n_1
\mathbf{1}$ is the dual poset of
$\mathbb{P}_0$.\\
\st Let $\{e_1 = (1, 0, \cdots , 0), ~e_2 = (0, 1, \cdots , 0),
\cdots , ~e_n = (0, \cdots , 0, 1)\}$ be the standard basis of
$\mathbb{F}_q^n$. Let $g\in G = Aut(\mathbb{F}_q^n,
w_{\mathbb{P}_0})$. Then, for each $1\leq s\leq t$, and each $i\in
n_s \mathbf{1}$,
\begin{equation}
  g(e_i) = a_i e_{\rho_s(i)} + \sum_{l=1}^{n_1+\cdots
  +n_{s-1}}b_{li}e_l,
\end{equation}
where $a_i\in \mathbb{F}_q^{\times}$, $b_{li}\in \mathbb{F}_q$,
$\rho_s(i)\in n_s \mathbf{1}$, and the assignment $i\mapsto
\rho_s(i)$ is a permutation on $n_s \mathbf{1}$. Conversely, if $g$
is the linear map given by (10), for $1\leq s\leq t$, $i\in n_s
\mathbf{1}$, then $g\in G = Aut(\mathbb{F}_q^n, w_{\mathbb{P}_0})$.
Let $x\cdot y = \sum_{i=1}^n x_i y_i$, for $x = (x_1, \cdots ,
x_n)$, $y = (y_1, \cdots , y_n)\in \mathbb{F}_q^n$. Then $\langle
x,y\rangle = \lambda(x\cdot y)$ is an inner product on
$\mathbb{F}_q^n$, where $\lambda$ is a fixed nontrivial additive
character on $\mathbb{F}_q$. Now, we will show the existence of a
map $\iota : G = Aut(\mathbb{F}_q^n, w_{\mathbb{P}_0})\rightarrow
\check{G} = Aut(\mathbb{F}_q^n, w_{\check{\mathbb{P}}_0})$ such that
$\langle gx,y\rangle = \langle x,\iota(g)y\rangle$, for all $g\in
G$, $x,y\in \mathbb{F}_q^n$. For this, it is enough to show that
there is a map $\iota : G\rightarrow \check{G}$ satisfying
\begin{equation}
  ge_i\cdot e_k = e_i\cdot \iota(g)e_k, ~for ~all ~g\in G, ~i,k = 1,
  \cdots , n.
\end{equation}
Let $g\in G$ be the map given by (10), for each $1\leq s\leq t$, and
each $i\in n_s \mathbf{1}$. If (11) is to be satisfied, we must have :\\
for $1\leq j\leq t$, $k\in n_j \mathbf{1}$,
\begin{align*}
  \iota(g)e_k\cdot e_i &= ge_i\cdot e_k\\
    &= (a_ie_{\rho_s(i)}+\sum_{l=1}^{n_1+ \cdots
    +n_{s-1}}b_{li}e_l)\cdot e_k\\
    &=
        \begin{cases}
        0,&\text{if} ~s<j,\\
        a_i,&\text{if} ~s=j ~and ~k=\rho_s(i),\\
        0,&\text{if} ~s=j ~and ~k\neq \rho_s(i),\\
        b_{ki},&\text{if} ~s>j.
        \end{cases}
\end{align*}
This yields that, for $1\leq j\leq t$, $k\in n_j \mathbf{1}$,
\begin{equation}
  \iota(g)e_k = a_{\rho_j^{-1}(k)}e_{\rho_j^{-1}(k)}+\sum_{l=n_1+ \cdots +n_j+1}^{n_1+ \cdots+ n_t}b_{kl}e_l.
\end{equation}
Now, $\iota(g)$ given by (12) belongs to $\check{G} =
Aut(\mathbb{F}_q^n, w_{\check{\mathbb{P}}_0})$, so that a map $\iota
: G\rightarrow \check{G}$ is defined and (11) is satisfied. Thus the
results stated in Corollary 10 hold true. In particular, $H(n_1,
n_2, \cdots ,n_t,q)$ and $H(n_t, n_{t-1}, \cdots ,n_1,q)$ are dual
to each other.


\begin{rem}\label{remark2}
(1) As we have seen in the above, the weak Hamming scheme $H(n_1,
\cdots ,n_t,q)$ is the associated scheme $\mathfrak{X}_G$ when $G =
Aut(\mathbb{F}_q^n, w_{\mathbb{P}_0}) ~(n = n_1+ \cdots +n_t)$ acts
on $\mathbb{F}_q^n$, with $\mathbb{P}_0 = n_1 \mathbf{1} \bigoplus
\cdots \bigoplus n_t \mathbf{1}$ the weak order poset. The weak
order poset is unique in many respects. Indeed, the following has
been shown. Let $\mathbb{P}$ be a poset on $[n]$. Then (1)
$\mathbb{P}$ is a weak order poset on $[n] \Leftrightarrow$ (2)
$(\mathbb{P},\check{\mathbb{P}})$is a weak dual MacWilliams pair
$(wdMp) \Leftrightarrow$ (3) The group $Aut(\mathbb{F}_q^n,
w_{\mathbb{P}})$ acts transitively on each $\mathbb{P}$-sphere
$S_{\mathbb{P}}(i) = \{x\in \mathbb{F}_q^n|w_{\mathbb{P}}(x) = i\}
~(0\leq i\leq n) \Leftrightarrow$ (4) $\mathfrak{X} =
(\mathbb{F}_q^n, ~{\{R_{i}\}}_{i=0}^{n})$,  with $(x,y)\in R_i
\Leftrightarrow x-y\in S_{\mathbb{P}}(i) ~(0\leq i\leq n)$, is an
association scheme. Here $\check{\mathbb{P}}$ is the dual poset of
$\mathbb{P}$, a pair of posets $(\mathbb{P},\check{\mathbb{P}})$ on
$[n]$ is called a $wdMp$ if the $\mathbb{P}$-weight distribution of
$C$ uniquely determines $\check{\mathbb{P}}$-weight distribution of
$C^{\perp}$, for every linear code $C\subseteq \mathbb{F}_q^n$. For
details about these, one is referred to [9-11,13,14].\\
(2) The structure of $Aut(\mathbb{F}_q^n, w_{\mathbb{P}_0})$ was
explicitly determined in \cite{DSKim:2007}. So the map $\iota : G =
Aut(\mathbb{F}_q^n, w_{\mathbb{P}_0})\rightarrow  \check{G}=
Aut(\mathbb{F}_q^n, w_{\check{\mathbb{P}}_0})$ would have been more
explicitly determined, just as we did in the example (g) of Section
4.
\end{rem}

\end{document}